\documentclass{article}\begin{document}\centerline{\bf The First Derivative of Ramanujans Cubic Continued Fraction}\vskip .10in

\centerline{\bf Nikos Bagis}

\centerline{Department of Informatics}

\centerline{Aristotele University of Thessaloniki Greece}
\centerline{nikosbagis@hotmail.gr}

\begin{quote}
\begin{abstract}
We give the complete evaluation of the first derivative of the Ramanujans cubic continued fraction using Elliptic functions. The Elliptic functions are easy to handle and give the results in terms of Gamma functions and radicals from tables.     
\end{abstract}

\bf keywords \rm{Ramanujan's Cubic Fraction; Jacobian Elliptic Functions; Continued Fractions; Derivative}

\end{quote}

\section{Introduction}
\label{intro}
The Ramanujan's Cubic Continued Fraction is (see [3], [7], [8], [9], [11]).
\begin{equation}
V(q):=\frac{q^{1/3}}{1+}\frac{q+q^2}{1+}\frac{q^2+q^4}{1+}\frac{q^3+q^6}{1+}\ldots
\end{equation}
Our main result is the evaluation of the first derivative of  Ramanujan's cubic fraction. For this, we follow a different way from previous works and use the theory of Elliptic functions. Our method consists to find the complete polynomial equation of the cubic fraction which is a solvable, in radicals, quartic equation, in terms only of the inverse elliptic nome $k_r$,  Using the derivative of $k_r$ which we evaluate in Section 2 of this article, we find the desired formula of the first derivative. For beginning we give some definitions first.\\
Let           
\begin{equation}
\left(a;q\right)_k=\prod^{k-1}_{n=0}(1-aq^n)
\end{equation}
Then we define
\begin{equation}
f(-q)=(q;q)_\infty 
\end{equation}
and
\begin{equation}
\Phi(-q)=(-q;q)_\infty
\end{equation}
Also let
\begin{equation}
K(x)=\int^{\pi/2}_{0} \frac{1}{\sqrt{1-x^2\sin^2(t)}}dt
\end{equation}
be the elliptic integral of the first kind.\\
We denote 
\begin{equation}
\theta_4(u,q)=\sum^{\infty}_{n=-\infty}(-1)^nq^{n^2}e^{2nui}
\end{equation}
the Elliptic Theta function of the 4th-kind. Also hold the following relations (see [16]):
\begin{equation}
\prod^{\infty}_{n=1}(1-q^{2n})^6=\frac{2kk'K(k)^3}{\pi^3q^{1/2}}
\end{equation}
\begin{equation}
q^{1/3}\prod^{\infty}_{n=1}(1+q^n)^8=2^{-4/3}\left(\frac{k}{1-k^2}\right)^{2/3} 
\end{equation}
and
\begin{equation}
f(-q)^8=\prod^{\infty}_{n=1}(1-q^n)^8=\frac{2^{8/3}}{\pi^4}q^{-1/3}k^{2/3}(k')^{8/3}K(k)^4 
\end{equation}
The variable $k$ is defined from the equation  
\begin{equation}
\frac{K(k')}{K(k)}=\sqrt{r}
\end{equation}
where $r$ is positive , $q=e^{-\pi \sqrt{r}}$ and $k'=\sqrt{1-k^2}$. Note also that whenever $r$ is positive rational, the $k=k_r$ are algebraic numbers. 
\section{The Derivative $\left\{r,k\right\}$} 
\textbf{Lemma 1.}\\
If $\left|t\right|<\pi a/2$ and $q=e^{-\pi a}$ then
\begin{equation}
\sum^{\infty}_{n=1}\frac{\cosh(2tn)}{n\sinh(\pi a n)}=\log(f(-q^2))-\log\left(\theta_4(it,e^{-a\pi})\right)
\end{equation}
\textbf{Proof.}\\
From the Jacobi Triple Product Identity (see [4]) we have 
\begin{equation}
\theta_4(z,q)=\prod^{\infty}_{n=0}(1-q^{2n+2})(1-q^{2n-1}e^{2iz})(1-q^{2n-1}e^{-2iz})
\end{equation}
By taking the logarithm of both sides and expanding the logarithm of the individual terms in a power series it is simple to show (11) from (12).
\[
\]
\textbf{Lemma 2.}\\
Let $q=e^{-\pi \sqrt{r}}$ with $r$ real positive 
\begin{equation}
\phi(x)=2\frac{d}{dx}\left(\frac{\partial}{\partial t}\log\left(\vartheta_4\left(\frac{it\pi}{2},e^{-2\pi x}\right)\right)_{t=x}\right)
\end{equation}
then
\begin{equation}
\frac{d(\sqrt{r})}{dk}=\frac{K^{(1)}(k)}{\phi\left(\frac{K(\sqrt{1-k^2})}{K(k)}\right)}=
\frac{K^{(1)}(k)}{\phi\left(\frac{K(k')}{K(k)}\right)}
\end{equation}
Where $K^{(1)}(k)$ is the first derivative of $K$.\\
\textbf{Proof.}\\
From Lemma 1 we have 
$$2\frac{\partial}{\partial t}\log\left(\vartheta_4\left(\frac{it\pi}{2},e^{-2\pi x}\right)\right)_{t=x}=-\pi\sum^{\infty}_{n=1}\frac{1}{\cosh\left(n\pi x\right)}=\frac{\pi}{2}-K(k_x)$$
then
$$\sqrt{x(k_2)}-\sqrt{x(k_1)}=-\int^{k_2}_{k_1}\frac{K^{(1)}(k)}{\phi\left(\frac{K\left(\sqrt{1-k^2}\right)}{K(k)}\right)}dk$$ 
Differentiating the above relation with respect to $k$ we get the result.
\[
\]
\textbf{Lemma 3.}\\
Set $q=e^{-\pi \sqrt{r}}$ and 
$$\left\{r,k\right\}:=\frac{dr}{dk}=2\frac{K(k')K^{(1)}(k)}{K(k)\phi\left(\frac{K(k')}{K(k)}\right)}$$
Then
\begin{equation}
\left\{r,k\right\}=\frac{\pi\sqrt{r}}{K^2(k_r)k_rk'^2_r}
\end{equation}
\textbf{Proof.}\\
From (9) taking the logarithmic derivative with respect to $k$ and using Lemma 2 we get:
\begin{equation}
\pi\left\{r,k\right\}\left(1-24\sum^{\infty}_{n=1}\frac{nq^n}{1-q^n}\right)=\left(\frac{1-5k^2}{(k-k^3)}+\frac{6K^{(1)}}{K}\right)\frac{4K'}{K}
\end{equation}
But it is known that
\begin{equation}
\sum^{\infty}_{n=1}\frac{nq^n}{1-q^n}=\frac{1}{24}+\frac{K}{6\pi^2}((5-k^2)K-6E)
\end{equation}
Hence
\begin{equation}
\left\{r,k\right\}=\frac{\pi K'}{K^2}\frac{\frac{1-5k^2}{k-k^3}+\frac{6K^{(1)}}{K}}{(k^2-5)K+6E}
\end{equation}
Also $$a(r)=\frac{\pi}{4K^2}+\sqrt{r}-\frac{E\sqrt{r}}{K} ,$$
where $a(r)$ is the elliptic alpha function. Using the above relations we get the result.\\
\textbf{Note.}\\
1) The first derivative of $K$ is $$K^{(1)}=\frac{E}{k_r\cdot k'^2_r}-\frac{K}{k_r}$$
where $k=k_r$ and $k'=k'_r=\sqrt{1-k^2_r}$.\\ 
2) In the same way we can find form the relation
\begin{equation}
k_{4r}=\frac{1-k'_r}{1+k'_r}
\end{equation} 
the 2-degree modular equation of the derivative.\\ 
Noting first that (the proof is easy)
\begin{equation}
\left\{r,k'_{r}\right\}=\frac{k'_r}{k_{r}}\left\{r,k_r\right\}
\end{equation}
we have
\begin{equation}
\left\{r,k_{4r}\right\}=\frac{k'_r(1+k'_r)^2}{2k_r}\left\{r,k_r\right\}
\end{equation}
\section{The Ramanujan's Cubic Continued Fraction}
Let 
\begin{equation} V(q):=\frac{q^{1/3}}{1+}\frac{q+q^2}{1+}\frac{q^2+q^4}{1+}\frac{q^3+q^6}{1+}\ldots
\end{equation}
is the Ramanujan's cubic continued fraction, then holds
\[
\]
\textbf{Lemma 4.}
\begin{equation}
V(q)=\frac{2^{-1/3}(k_{9r})^{1/4}(k'_{r})^{1/6}}{(k_r)^{1/12}(k'_{9r})^{1/2}}
\end{equation}
where the $k_{9r}$ are given by (see [7]):
\begin{equation}
\sqrt{k_rk_{9r}}+\sqrt{k'_rk'_{9r}}=1
\end{equation}
\textbf{Proof.}\\
The proof can be found in [18].
\[
\]
\textbf{Lemma 5.}\\
If $$G(x)=\frac{x}{\sqrt{2\sqrt{x}-3x+2x^{3/2}-2\sqrt{x}\sqrt{1-3\sqrt{x}+4x-3x^{3/2}+x^2}}}$$ and 
\begin{equation}
k_r=G(w)
\end{equation}
then
$$k_{9r}=\frac{w}{k_r}$$
and
$$k'_{9r}=\frac{(1-\sqrt{w})^2}{k'_r}$$
\textbf{Proof.}\\
See [18]. 
\[
\]
\textbf{Theorem 1.}\\
Set $T=\sqrt{1-8V^3(q)}$ then
\begin{equation}
(k_r)^2=\frac{(1-T)(3+T)^3}{(1+T)(3-T)^3}
\end{equation}
\textbf{Proof.}\\
See [18].
\[
\]
Equation (26) is a solvable quartic equation with respect to $T$.\\ 
An example of evaluation is
\begin{equation}
V(e^{-\pi})=\frac{1}{2}\left(-2-\sqrt{3}+\sqrt{3(3+2\sqrt{3})}\right)
\end{equation} 
\[
\]
\textbf{Main Theorem.}\\
Let $q=e^{-\pi\sqrt{r}}$, then
\begin{equation}
V'(q)=\frac{dV(q)}{dq}=\frac{-2\sqrt{r}}{q\pi}\frac{dV}{dr}=\frac{4K^2(k_r)k'^2_r(V(q)+V^4(q))}{3q\pi^2\sqrt{r}\sqrt{1-8V^3(q)}}
\end{equation}  
\textbf{Proof.}\\ Derivate (26) with respect to $r$ then 
\begin{equation}
\sqrt{\frac{2k_r}{\left\{k,r\right\}}}=\frac{4T(3+T)}{(3-T)^2(1+T)}\sqrt{\frac{dT}{dr}}
\end{equation}
or
\begin{equation}
T_r=\frac{dT}{dr}=\frac{1}{8k_r\left\{r,k\right\}}\frac{(9-T^2)(1-T^2)}{T^2}
\end{equation}
Using the relation $T=\sqrt{1-8V(q)^3}$, we get
$$
\frac{dV(q)}{dr}=-\frac{2}{3}\frac{V(q)+V^4(q)}{k_r \left\{r,k\right\} \sqrt{1-8V^3(q)}}
\eqno{(a)}$$
which is the result.\\ Hence the problem of finding $V(q)$ and $V'(q)$ is completely solvable in radicals when we know $k_r$ and $K(k_r)$ (see [12]), $r\in\bf Q\rm$, $r>0$.
\[
\]
We often use the notations $V[r]:=V(e^{-\pi\sqrt{r}})$,  $T[r]:=T(e^{-\pi\sqrt{r}})=t$.
\[
\] 
\textbf{Proposition 1.}
\begin{equation}
V[4r]=\frac{1-T[r]}{4V[r]}
\end{equation}
\textbf{Proof.}\\
See [9].
\[
\]
\textbf{Proposition 2.}\\
Set $T'[4r]=u$, $T'[r]=\nu$, then 
\begin{equation}
\frac{u}{\nu}=\frac{(1-t)(3+t)}{8\sqrt{t}(1+t)^{5/3}(3-t)^{1/2}}  
\end{equation}
\textbf{Proof.}\\
From (19), (20), (21) and (28) we get  
\begin{equation}
V'[4r]=\frac{-2\left\{k,r\right\}}{3\frac{1-k'}{1+k'}\frac{k'(1+k')^2}{2k}}\frac{V[4r]+V[4r]^4}{T[r]}
\end{equation}
If we use the duplication formula (31) we get the result.
\[
\]   
\textbf{Evaluations}.\\
1) We can calculate now easy the values of $V'(q)$  from (28) using (26). An example of evaluation is $$k_1=\frac{1}{\sqrt{2}}$$, $$E(k_1)=\frac{4\pi^{3/2}}{\Gamma(-1/4)^2}+\frac{\Gamma(3/4)^2}{2\sqrt{\pi}}$$ and $$K(k_1)=\frac{8\pi^{3/2}}{\Gamma(-1/4)^2}$$
When $r=1$ we get $$\left\{r,k\right\}=\frac{8\sqrt{2}\Gamma(3/4)^4}{\pi^2}$$ Hence 
\begin{equation}
V'(e^{-\pi})=-\frac{64 \left(-26-15 \sqrt{3}+10 \sqrt{3+2 \sqrt{3}}+6 \sqrt{9+6 \sqrt{3}}\right)}{\sqrt{45+26 \sqrt{3}-18 \sqrt{3+2 \sqrt{3}}-10 \sqrt{9+6 \sqrt{3}}}}\frac{e^{\pi } \pi}{\Gamma\left(-\frac{1}{4}\right)^4}
\end{equation} 
2) It is $$T_1=T(e^{-\pi\sqrt{3}})=-39+22\sqrt{3}-\frac{2\cdot 6^{2/3}(-123+71\sqrt{3})}{\left(-4725+2728\sqrt{3}-\sqrt{4053-2340\sqrt{3}}\right)^{1/3}}+
$$
$$
+2\cdot6^{1/3}\left(-4725+2728\sqrt{3}-\sqrt{4053-2340\sqrt{3}}\right)^{1/3}$$
and $V_1=V(e^{-\pi\sqrt{3}})=\frac{1}{2}\sqrt[3]{1-T_1^2}$.\\ From tables and (15) it is:
$$
\left\{3,k_3\right\}=\frac{192 \sqrt{2} \left(-1+\sqrt{3}\right) \pi ^2}{\Gamma\left(\frac{1}{6}\right)^2 \Gamma\left(\frac{1}{3}\right)^2}$$
We find the value of $V'(e^{-\pi\sqrt{3}})$ in terms of Gamma function and algebraic numbers. 
$$V'(e^{-\pi\sqrt{3}})=\frac{4\sqrt{3}e^{\pi\sqrt{3}}}{3}\frac{V_1+V^4_1}{k_3 \left\{3,k_3\right\} \sqrt{1-8 V^3_1}}$$

\[
\]

\centerline{\bf References}\vskip .2in

\noindent

[1]: M.Abramowitz and I.A.Stegun. 'Handbook of Mathematical Functions'. Dover Publications, New York. 1972.

[2]: C. Adiga, T. Kim. 'On a Continued Fraction of Ramanujan'. Tamsui Oxford Journal of Mathematical Sciences 19(1) (2003) 55-56 Alethia University.

[3]: C. Adiga, T. Kim, M.S. Naika and H.S. Madhusudhan. 'On Ramanujan`s Cubic Continued Fraction and Explicit Evaluations of Theta-Functions'. arXiv:math/0502323v1 [math.NT] 15 Feb 2005.

[4]: G.E.Andrews. 'Number Theory'. Dover Publications, New York. 1994.

[5]: B.C.Berndt. 'Ramanujan`s Notebooks Part I'. Springer Verlag, New York (1985).

[6]: B.C.Berndt. 'Ramanujan`s Notebooks Part II'. Springer Verlag, New York (1989).

[7]: B.C.Berndt. 'Ramanujan`s Notebooks Part III'. Springer Verlag, New York (1991).

[8]: Bruce C. Berndt, Heng Huat Chan and Liang-Cheng Zhang. 'Ramanujan`s class invariants and cubic continued fraction'. Acta Arithmetica LXXIII.1 (1995).   

[9]: Heng Huat Chan. 'On Ramanujans Cubic Continued Fraction'. Acta Arithmetica. 73 (1995), 343-355.     

[10]: I.S. Gradshteyn and I.M. Ryzhik. 'Table of Integrals, Series and Products'. Academic Press (1980).

[11]: Megadahalli Sidda Naika Mahadeva Naika, Mugur Chin

[12]: Habib Muzaffar and Kenneth S. Williams. 'Evaluation of Complete Elliptic Integrals of The First Kind and Singular Moduli'. Taiwanese Journal of Mathematics. Vol. 10, No. 6, pp, 1633-1660, Dec 2006.   

[13]: L. Lorentzen and H. Waadeland. 'Continued Fractions with Applications'. Elsevier Science Publishers B.V., North Holland (1992).  

[14]: S.H.Son. 'Some integrals of theta functions in Ramanujan's lost notebook'. Proc. Canad. No. Thy Assoc. No.5 (R.Gupta and K.S.Williams, eds.), Amer. Math. Soc., Providence.

[15]: H.S. Wall. 'Analytic Theory of Continued Fractions'. Chelsea Publishing Company, Bronx, N.Y. 1948.  

[16]:E.T.Whittaker and G.N.Watson. 'A course on Modern Analysis'. Cambridge U.P. (1927)

[17]:I.J. Zucker. 'The summation of series of hyperbolic  functions'. SIAM J. Math. Ana.10.192 (1979)

[18]: Nikos Bagis. 'The complete evaluation of Rogers Ramanujan and other continued fractions with elliptic functions'. arXiv:1008.1304v1  

\end{document}